\documentclass[12pt]{article}
\usepackage{standalone}
\usepackage{amsmath, amsthm, amssymb, color}
\usepackage[colorlinks=true,linkcolor=blue,urlcolor=blue]{hyperref}
\usepackage{graphicx}
\usepackage{caption}
\usepackage{mathtools}
\usepackage{enumerate}
\usepackage{verbatim}
\usepackage{tikz,tikz-cd,tikz-3dplot}
\usepackage{amssymb}
\usetikzlibrary{matrix}
\usetikzlibrary{arrows}
\usepackage{algorithm}
\usepackage[noend]{algpseudocode}
\usepackage{caption}
\usepackage[normalem]{ulem}
\usepackage{subcaption}
\usepackage{multicol}
\oddsidemargin \evensidemargin
\usepackage{makecell}
\usepackage{array}
\usepackage{enumitem}

\newtheorem{theorem}{Theorem}[section]

\newtheorem{proposition}[theorem]{Proposition}
\newtheorem{lemma}[theorem]{Lemma}
\newtheorem{corollary}[theorem]{Corollary}

\theoremstyle{definition}

\newtheorem{remark}[theorem]{Remark}

\newtheorem{examplex}[theorem]{Example}
\newenvironment{example}
{\pushQED{\qed}\examplex}
{\popQED\endexamplex}

\newcommand{\bR}{\mathbb R}
\newcommand{\bC}{\mathbb C}
\newcommand{\bZ}{\mathbb Z}

\newcommand{\bP}{\mathbb P}
\newcommand{\bS}{\mathbb S}

\newcommand{\cB}{\mathcal B}

\newcommand{\cO}{\mathcal O}
\newcommand{\cT}{\mathcal T}
\newcommand{\cV}{\mathcal V}

\renewcommand{\phi}{\varphi}

\newcommand{\gen}[1]{\langle {#1} \rangle}

\newcommand{\Gr}{\text{Gr}}

\newcommand{\Sym}{\text{Sym}}

\newcommand{\SL}{\text{SL}}
\newcommand{\vect}[2]{\begin{bmatrix} #1 \\ #2 \end{bmatrix}}
\newcommand{\Seg}{\text{Seg}}

\usepackage{xcolor}

\usepackage[breakable,skins]{tcolorbox}
\newtcolorbox{tbox}[1][]{%
	breakable,
	enhanced,
	colframe=black,
	coltitle=white,
	#1
}

\title{\bf The Segre Determinant}
\author{Elizabeth Pratt\footnote{Supported by NSF GRFP no. 2023358166 }}
\date{}
\begin{document}
	\maketitle
	
	\begin{abstract}
		The Segre determinant is a polynomial which encodes the condition for points to lie on a bilinear hypersurface in the product of projective spaces. We study Segre determinants and compute them in various coordinate systems. We show that the Segre determinant represents the Chow-Lam form of a generic torus orbit in the Grassmannian. These Chow-Lam forms were introduced as a generalization of Chow forms for projective varieties, and enjoy many similar properties. We also present applications to algebraic vision and to Chow quotients of Grassmannians.
	\end{abstract}
	
	\section{Introduction}
	
	Fix vector spaces $V$ and $W$ over $\bC$ of dimensions $k$ and $\ell.$  Consider $r$ points $A_1 \times B_1, \, \ldots, \,  A_r \times B_r$ in $\bP(V) \times \bP(W)$. The \emph{Segre matrix} of this point configuration is the $k\ell \times r$ matrix
	\begin{equation}\label{eq:segmatrix}
		\begin{bmatrix}
			\vdots & & \vdots \\
			A_1 \otimes B_1 & \hdots & A_r \otimes B_r \\
			\vdots & & \vdots
		\end{bmatrix} .
	\end{equation}
	It is a flattening of the $k \times \ell \times r$ tensor with slices $A_i \otimes B_i.$ When $r = k\ell,$ the Segre matrix is square and we call its determinant the \emph{Segre determinant}. This determinant vanishes whenever the $k\ell$ points lie on a hyperplane section of $\bP(V) \times \bP(W)$ in the Segre embedding, hence the name. Indeed, the equation of such a hyperplane section is in the left kernel of the Segre matrix. We will use $\Seg_{k,\ell}$ to denote the Segre determinant as a polynomial in the coordinates of the points $A_i$ and $B_i.$
	
	The aim of this paper is to introduce the Segre determinant and give various applications. After establishing its basis properties (Section \ref{sec:computation}), we show how it appears in various guises: in algebraic vision (Section \ref{sec:vision}), as a Chow-Lam form (Sections \ref{sec:chowlam}--\ref{sec:torus}), and in the theory of Chow quotients of Grassmannians (Section \ref{sec:chowquotient}).  
	
	The Segre determinant is a generalization of the classical determinant for $k = 1.$ The special case $k=\ell= 2$ encodes the condition for two configurations of four ordered points in $\bP^1$ to be the same up to automorphism; see Example \ref{eg:22} below. The case $k=\ell = 3$ appears in algebraic vision as a necessary condition for two configurations of nine ordered points in $\bP^2$ to be linear projections of a common configuration of ordered points in $\bP^3$ \cite{Agarwal}.
	
	\begin{example}\label{eg:22}
		Consider four points $\vect{a_{1,1}}{a_{2,1}} \times \vect{b_{1,1}}{b_{2,1}}, \, \dots, \, \vect{a_{1,4}}{a_{2,4}} \times \vect{b_{1,4}}{b_{2,4}}$ in the surface $\bP^1 \times \bP^1.$ The Segre matrix looks like
		\begin{equation}
			\begin{bmatrix}
				a_{1,1}b_{1,1} & a_{1,2}b_{1,2} & a_{1,3}b_{1,3} & a_{1,4}b_{1,4} \\
				a_{1,1}b_{2,1} & a_{1,2}b_{2,2} & a_{1,3}b_{2,3} & a_{1,4}b_{2,4} \\
				a_{2,1}b_{1,1} & a_{2,2}b_{1,2} & a_{2,3}b_{1,3} & a_{2,4}b_{1,4} \\
				a_{2,1}b_{2,1} & a_{2,2}b_{2,2} & a_{2,3}b_{2,3} & a_{2,4}b_{2,4}
			\end{bmatrix}.
		\end{equation}
		This expression vanishes whenever the four points in $\bP^1 \times \bP^1$ lie on a common $(1,1)$-curve. We collect our eight vectors into two matrices
		\[A \ := \begin{bmatrix}
			a_{1,1} & a_{1,2} & a_{1,3} & a_{1,4} \\
			a_{2,1} & a_{2,2} & a_{2,3} & a_{2,4}
		\end{bmatrix}, \qquad B \ := \begin{bmatrix} b_{1,1} & b_{1,2} & b_{1,3} & b_{1,4} \\
			b_{2,1} & b_{2,2} & b_{2,3} & b_{2,4}
		\end{bmatrix}.\]
		Let $[ij]$ denote the six maximal minors of $A$ and $\gen{ij}$ denote the six maximal minors of $B.$ Then we compute
		\[\text{Seg}_{2,2} = [12][34]\gen{13}\gen{24}-[13][24]\gen{12}\gen{34}.\]
		It vanishes whenever the cross-ratio $\frac{[12][34]}{[13][24]}$ of the four points $A_i$ in $\bP^1$ is equal to the cross-ratio of the four points $B_i$ in $\bP^1.$ The cross-ratio is the fundamental invariant for $\text{PGL}_2$ acting on $(\bP^1)^4.$
	\end{example}
	
	The first two sections focus on computations and examples. In Section \ref{sec:computation} we prove that the Segre determinant may be written in terms of  the maximal minors of the $k \times k\ell$ matrix $A$ and the $\ell \times k\ell$ matrix $B,$ generalizing the computation in Example \ref{eg:22}. These minors generate the invariant rings of $\SL_k$ and $\SL_\ell$ acting on $(\bP^{k-1})^n$ and $(\bP^{\ell-1})^n,$ respectively \cite[Theorem 3.2.1]{AIT}. We also discuss how to compute such expressions, which is difficult in practice. 
	
	This flavor of question has a long history in classical algebraic geometry, and is related to the problem of synthetic (ruler-and-compass) constructions. The analogous problem for the Veronese embedding is to figure out when $\binom{n+k-1}{k}$ points in $\bP^{n-1}$ lie on a hypersurface of degree $k.$ Pascal's theorem for $k=2$ and $n=3$ is the earliest such result, giving a condition for six points in the plane to lie on a conic. The \emph{Bruxelles problem}, posed in 1825 by l’Acad\'emie de Bruxelles, asks for a synthetic construction to determine when ten points in $\bP^3$ lie on a quadric surface ($k=2, \, n=4$). This is a significant leap in difficulty, and a construction was only recently obtained by Traves \cite{Traves}. 
	
	Such a property of a point configuration is $SL_k$-invariant, so we may ask how to write it in terms of minors of $A$ and $B$, perhaps thus gleaning some geometric insight. Indeed, in 1927 Turnbull and Young gave an algebraic expression in terms of the $4 \times 4$ minors of the $4 \times 10$ matrix parameterizing the points \cite{TurnbullYoung}. In this spirit, we present the analogous computation for eight points in $\bP^1 \times \bP^3$ and nine points in $\bP^2 \times \bP^2$ in Equations \eqref{eq:seg24} and \eqref{eq:seg33}, respectively.
	
	Section \ref{sec:vision} gives an application of the Segre determinant $\Seg_{3,3}$ to algebraic vision. In this setting, the Segre determinant is a necessary condition for two distinct configurations of nine ordered points in $\bP^2$ to be linear projections of a common set of nine ordered points in $\bP^3.$ Our computation answers a question of Rekha Thomas about writing this condition in terms of $SL_3$-invariants. This section is completely independent from Sections \ref{sec:chowlam}--\ref{sec:chowquotient}.
	
	Sections \ref{sec:chowlam} and \ref{sec:torus} give a geometric interpretation of the Segre determinant for general $k$ and $\ell$ as the Chow-Lam form of a torus orbit in $\Gr(k,k\ell).$ The Chow-Lam form was introduced in \cite{ChowLam} as a generalization of the classical Chow form in algebraic geometry. It associates to a subvariety of $\Gr(k,n)$ a single polynomial which in most cases completely encodes the variety \cite{PrattRanestad}. These forms first arose in calculations of scattering amplitudes in high-energy physics. Our main theorem in Section \ref{sec:torus} is the following.
	
	\begin{theorem}[Segre Determinant]\label{thm:intro}
		Suppose $k, l \geq 2$ and let $n = k\ell.$ Fix a point $A$ in ${\rm Gr}(k,n)$ with non-zero Pl\"ucker coordinates. Then the Chow-Lam form of the torus orbit of $A$ in primal Pl\"ucker coordinates $B$ on ${\rm Gr}(n-l, n)$ equals the Segre determinant ${\rm Seg}_{k,\ell}(A,B).$
	\end{theorem}
	
	In Section \ref{sec:chowquotient} we define the \emph{Segre coefficient variety}. It parameterizes Segre determinants as polynomials in the $B$ variables as $A$ varies. We show in Theorem \ref{thm:independent} that every monomial of the form $[I_1] \cdots [I_\ell]$ with $I_1 \cup \ldots \cup I_\ell = [k\ell]$ appears in the linear span of the coefficients of the Segre polynomial. This results in Corollary \ref{cor:git}, which states that the Segre coefficient variety for $k=2$ is isomorphic to the GIT quotient $(\bP^1)^{2 \ell} / \!/ \text{SL}(2).$ Thus the Segre determinant uniquely determines torus orbit closures of generic points in $\Gr(2,n).$ However, this fails for $k=3.$ We give two torus orbit closures with the same Segre determinant in Example \ref{eg:degree2}.
	
	\subsection{Acknowledgements}
	I am grateful to Bernd Sturmfels for many helpful discussions and feedback on this work. I also thank Rekha Thomas and Timothy Duff for teaching me about algebraic vision, and Will Traves for useful discussions regarding synthetic constructions. Thanks also to the anonymous reviewer for their feedback.
	
	\section{Properties and Computation} \label{sec:computation}
	We now discuss the computation of the Segre determinant in various coordinate systems. As in the introduction, fix vector spaces $V$ and $W$ of dimensions $k$ and $\ell,$ and let $n := k \ell.$  
	
	Consider the polynomial rings $\bC[a_{ij}]_{1 \leq i \leq k, 1 \leq j \leq n}$ and $\bC[b_{ij}]_{1 \leq i \leq \ell, 1 \leq j \leq n}.$  As in the introduction, we may  collect the indeterminates $a_{ij}$ and $b_{ij}$ into a $k \times n$ matrix $A$ and an $\ell \times n$ matrix $B.$ Given a subset $I \in \binom{[n]}{k},$ we define the \emph{bracket} $[I]$ to be the determinant of the submatrix of $A$ with columns indexed by $I.$ Similarly, given $J \in \binom{[n]}{\ell},$ we define the bracket $\gen{J}$ to be the determinant of the submatrix of $B$ with columns $J.$
	These rings have actions of $\SL_k$ and $\SL_\ell,$ respectively, given by left multiplication of the matrices $A$ and $B.$ The following is sometimes called the \emph{First Fundamental Theorem of Invariant Theory}.
	\begin{theorem}[Theorem 3.2.1 of \cite{AIT}]\label{thm:fundinv}
		The brackets $[I]$ generate the invariant ring
		$\bC[a_{ij}]^{\SL_k}.$
	\end{theorem}	
	
	We denote the invariant ring $\bC[a_{ij}]^{\SL_k}$ by $\cB_{k,n}$ and call it the \emph{bracket algebra}. The Segre determinant lives in the tensor product  $\bC[a_{ij}] \otimes \bC[b_{ij}]$ of $\bC$-algebras. It is separately invariant under the $\SL_k$ and $\SL_\ell$ actions on $V \otimes W$, and thus may be written either in brackets $[I]$ and indeterminates $b_{ij},$ or in brackets $\gen{J}$ and indeterminates $a_{ij}.$ These can be realized concretely as two different block Laplace expansions of the Segre matrix. The following result states that $\Seg_{k,\ell}$ may be written \emph{simultaneously} in the two systems of brackets.
	
	\begin{proposition}\label{prop:brackets}
		The Segre determinant ${\rm Seg}_{k,\ell}$ is a polynomial of bi-degree $(l, k)$ in the brackets $[I]$ and $\gen{J}.$
	\end{proposition}
	The main ingredient in proving Proposition \ref{prop:brackets} is the following lemma. 
	
	\begin{lemma}\label{lem:invariant}
		Suppose that $V$ and $W$ are finite-dimensional representations of groups $G$ and $H,$ respectively. Then
		\begin{equation} \label{eq:intersection}
			V^G \otimes W^H  = (V \otimes W)^{G \times H}.
		\end{equation}
	\end{lemma}
	\begin{proof}
		The inclusion $\subseteq$ is immediate. For the other inclusion, choose a basis $\alpha_1, \, \ldots, \, \alpha_r$ for $V^G$ and add vectors $v_j$ to extend it to a basis for $V.$ Similarly, choose a basis $\beta_1, \, \ldots, \, \beta_s$ for $W^H$ and add vectors $w_j$ to extend it to a basis for $W.$ We observe that the right-hand side of (\ref{eq:intersection}) is contained in the vector space
		\[\left(V^G \otimes W\right) \cap \left(V \otimes W^G \right) \ \subseteq \ V \otimes W.\]
		An element $f$ of the intersection may be written uniquely in each basis as
		\begin{align}
			f &= \sum_{i \leq r,j \leq s,} c_{ij} \alpha_i \otimes \beta_j + \sum_{i \leq r,j > s}^k c_{ij} \alpha_i \otimes w_j \\
			&= \sum_{i \leq r,j \leq s} d_{ij} \alpha_i \otimes \beta_j + \sum_{i > r,j \leq s}^k d_{ij} v_i \otimes \beta_j.
		\end{align}
		Since the expressions are unique, we must have that $f = \sum_{i \leq r,j \leq s} c_{ij} \alpha_i \otimes \beta_j.$
	\end{proof}
	
	\begin{proof}[Proof of Proposition~\ref{prop:brackets}]
		The proof is mostly type-checking. For clarity, within this proof we use $\boxtimes$ when forming the representation of a direct product of groups. The Segre determinant lives in the coordinate ring $\Sym^\bullet(V \boxtimes W) \otimes \cdots \otimes \Sym^\bullet(V \boxtimes W)$ of the product $\bP(V \boxtimes W) \times \cdots \times \bP(V \boxtimes W).$ This ring has a multi-grading given by taking the grading in each tensor factor, and the Segre determinant lives in the component with grading $1^n = (1, \, \ldots, \, 1)$. Consider the map of $\SL_k \times \SL_\ell$-representations
		\begin{equation}
			\begin{split}
				\psi: (V \boxtimes W)^{\otimes n}  & \to V^{\otimes n} \boxtimes W^{\otimes n} \\
				(x_1 \boxtimes y_1) \otimes \ldots \otimes (x_n \boxtimes y_n) & \mapsto (x_1 \otimes \ldots \otimes x_n) \boxtimes (y_1 \otimes \ldots \otimes y_n).
			\end{split} 
		\end{equation}
		There is a map $\SL_k \times \SL_l \to \text{GL}(V \boxtimes W)$ given by $(A, B) \cdot x_i \otimes y_j = Ax_i \otimes By_j.$ The image lies in $\SL_{k\ell} = \SL(V \boxtimes W).$ Thus there is an inclusion of invariant vector spaces 
		\begin{equation}\label{eq:slinclution}
			\iota: \left( (V \boxtimes W)^{\otimes n} \right)^{\SL_{k\ell}} \hookrightarrow (V^{\otimes n} \boxtimes W^{\otimes n})^{\SL_k \times \SL_\ell}.
		\end{equation}
		Finally, by Lemma \ref{lem:invariant} we have that 
		\begin{equation}\label{eq:boxinvariant}
			(V^{\otimes n} \boxtimes W^{\otimes n})^{\SL_k \times \SL_\ell} = (V^{\otimes n})^{\SL_k} \boxtimes (W^{\otimes n})^{\SL_\ell}.
		\end{equation}
		By the First Fundamental Theorem of Invariant Theory, the right-hand side of \eqref{eq:boxinvariant} is the algebra generated by the brackets $[I]$ and $\gen{J}.$ By \eqref{eq:slinclution} the Segre determinant $\Seg_{k,\ell}$ lies in $(V^{\otimes n} \boxtimes W^{\otimes n})^{\SL_k \times \SL_\ell}$, and applying $\iota$ gives its expansion into brackets $[I]$ and $\gen{J}.$ Note that the total degree of $\Seg_{k,\ell}$ in each set of variables $a_{ij}$ and $b_{ij}$ is $n.$ Since $[I]$ has degree $k$ and $\gen{J}$ has degree $\ell$ when expanded into these variables, the bi-degree of the Segre determinant is $(n / k, n / \ell) = (\ell, k).$
	\end{proof}
	
	\begin{remark}
		The Segre determinant has further symmetries. It transforms equivariantly under the actions of $(\bC^*)^n$ on $\cB_{k,n}$ and $\cB_{\ell,n}$ obtained by scaling columns of matrix representatives. It is also invariant (up to sign) under permuting the points $A_i \times B_i$ and, when $k=\ell,$ exchanging the roles of $A_i$ and $B_i$. These and the degree considerations in Proposition \ref{prop:brackets} are enough to deduce, for example, that the polynomial $\Seg_{2,2}$ in Example \ref{eg:22} is a linear combination of $[12][34]\gen{12}\gen{34} + [13][24]\gen{13}\gen{24}$ and $ [12][34]\gen{13}\gen{24}-[13][24]\gen{12}\gen{34}.$ 
	\end{remark}
	
	The brackets in $\cB_{k,n}$ satisfy certain relations called the Pl\"ucker relations, coming from the relations between maximal minors of a $k \times n$ matrix. As the name implies, these are the same relations defining the Grassmannian in its Pl\"ucker embedding. Indeed, $\text{Proj} \cB_{k,n}$ equals $\Gr(k,n).$ It is a fact, sometimes called the \emph{Second Fundamental Theorem of Invariant Theory}, that the Pl\"ucker relations generate all relations between brackets \cite[Theorem 3.1.7]{AIT}
	
	The bracket algebra comes with a convenient basis in each graded component, which we explain as follows. A \emph{Young tableau} is the filling of a partition diagram with entries in $[n],$ allowing repeats. A \emph{semi-standard Young tableau} has the additional property that the numbers are non-decreasing within each row and strictly increasing within each column. We call a monomial $[I_1] \cdots [I_r]$ in $\cB_{k,n}$  \emph{standard} if the $k \times r$ rectangular Young tableau obtained by stacking $I_1, \, \ldots, \, I_r$ vertically and then transposing is a semi-standard Young tableau. For example, in $\cB_{2,4}$ the monomials $[12][12]$ and $[13][24]$ are standard, but the monomial $[14][23]$ is not. There are $21$ degree two monomials in total in $\cB_{2,4}$, $20$ of which are standard. This use of the term ``standard monomial" is consistent with the theory of Gr\"obner bases: the Pl\"ucker relations form a Gr\"obner basis for the ideal of relations in $\cB_{k,n}$ under a certain order called tableau order, and these are the standard monomials \cite[Theorem 3.1.7]{AIT}. 
	
	The semi-standard monomials of a fixed degree $r$ form a basis for the degree $r$ part of the bracket algebra \cite[Corollary 3.1.9]{AIT}. The \emph{straightening algorithm} applies the Pl\"ucker relations to put a polynomial in $\cB_{k,n}$ into the unique representation such that every monomial within it is standard. Note that the Segre determinant is an eigenvector of the torus $(\bC^*)^{k \ell}$ acting on either $\cB_{k,k \ell}$ or $\cB_{\ell, k \ell}$. Thus each bracket monomial appearing within it has the same multi-set of indices, namely $[n]$. From here on when we speak of the ``standard basis" we mean the basis of semi-standard Young tableaux for the multi-linear component, in which the multi-set of the indices is $[n]$.
	
	\begin{example}\label{eg:seg24}
		The Segre determinant $\Seg_{2,4}$ has bi-degree $(4,2)$ and a total of $22$ monomials in the standard basis. Its expansion in standard brackets equals
		\begin{equation}\label{eq:seg24}
			\begin{aligned}
				\Seg_{2,4} \ = \ \small
				- \, (\gen{1235}\gen{4678} &- \gen{1245}\gen{3678} + \gen{1257}\gen{3468})[13][24][56][78] \\
				- \, (\gen{1237}\gen{4568} &- \gen{1357}\gen{2468} + \gen{1345}\gen{2678})[12][34][56][78] \\
				- \, (\gen{1235}\gen{4678} &- \gen{1236}\gen{4578} + \gen{1356}\gen{2478})[12][34][57][68] \\
				\, + \, \gen{1234}\gen{5678}[15][26][37][48] &
				\, + \, (\gen{1234}\gen{5678} + \gen{1346}\gen{2578})[12][35][47][68] \\
				- \gen{1347}\gen{2568}[12][35][46][78] & \, - \, \gen{1345}\gen{2678}[12][36][47][58] \\
				\, + \, \gen{1256}\gen{3478}[13][24][57][68] & \, - \, \gen{1246}\gen{3578}[13][25][47][68] \\
				+ \, \gen{1245}\gen{3678}[13][26][47][58] & \, - \, \gen{1237}\gen{4568}[14][25][36][78] \\
				+ \, \gen{1236}\gen{4578}[14][25][37][68] & \, - \, \gen{1235}\gen{4678}[14][26][37][58] \\
				+ \, (\gen{1234}\gen{5678} &+ \gen{1247}\gen{3568})[13][25][46][78].  &  \hfill \qedhere
			\end{aligned}
		\end{equation}
	\end{example}
	
	One may in principle compute expressions such as Example \ref{eg:seg24} by intersecting the expressions in $a_{ij}$ and $b_{ij}$ with the subring of invariants in $\bC[a_{ij}] \otimes \bC[b_{ij}]$. This intersection is typically computed via elimination algorithms which use Gr\"obner bases, such as that in \cite[Section 15.10.4]{Eisenbud}. However, these are unlikely to terminate. It is far more efficient to leverage the standard basis and use linear algebra. For instance, to do the above computation we performed the following steps.
	\begin{enumerate}
		\item Do a block Laplace expansion of the Segre matrix \eqref{eq:segmatrix} in the $[I]$ brackets and then straighten to the fourteen standard ones.
		\item For each standard monomial $T$ in the $[I]$ brackets:
		\begin{enumerate}[label=\roman*.]
			\item Define $f_T$ to be the coefficient of $T$ in the variables $b_{ij}$.
			\item Write $f_T = c_1\gen{1234}\gen{5678} + \, \ldots \, + c_{14}\gen{1347}\gen{2568}$, where the $c_i$ are unknown.
			\item Choose a random degree eight monomial $m$ in the $b_{ij}$ appearing in $f_T$. Its incidence vector with the brackets $\gen{J_1}\gen{J_2}$ imposes a linear constraint on the $c_i$.
			\item Repeat until the coefficients $c_i$ are determined. 
		\end{enumerate}
	\end{enumerate}
	The final step is linear algebra in a $14$-dimensional vector space, so it runs very quickly.
	
	\section{An Application to Algebraic Vision} \label{sec:vision}
	Computer vision is the study of ``cameras," namely linear projections from $\bP_{\bR}^3$ to $\bP_{\bR}^2,$ and how a computer gains information from them. A typical computer vision problem is to reconstruct an object in $3$-space from a set of $2$-dimensional snapshots. In algebraic vision, the object one is taking a picture of is an algebraic variety in $\bP_{\bR}^3,$ such as a curve or a collection of points. The survey \cite{Kileel} provides an overview of this research area. For the rest of this section we work over the field $\bR.$
	
	Fix two configurations $A$ and $B$ of eight ordered points in $\bP^2.$ One natural question in computer vision is: when are they linear projections of a common configuration in $\bP^3$? To answer this, consider the $9 \times 8$ Segre matrix with columns $A_1 \otimes B_1, \, \ldots, \, A_8 \otimes B_8$. This matrix has a one-dimensional kernel, provided that the point configurations are sufficiently generic (see \cite{LonguetHiggins} for the exact conditions). The kernel of the Segre matrix is a $1 \times 9$ vector. However, it may also be viewed as a $3 \times 3$ matrix, which we denote $F$. Then a necessary and sufficient condition for $A$ and $B$ to have a common recovery is that $F$ has rank two \cite{LonguetHiggins}. In algebraic vision $F$ is called the \emph{fundamental matrix}.
	
	If one has instead two configurations of \emph{nine} points, there is an extra condition: namely, the $9 \times 9$ Segre matrix must have a kernel. In Equation \eqref{eq:seg33} we express the Segre determinant as a polynomial in brackets $[I]$ and $\gen{J}.$ This answers a question of Rekha Thomas about how to write the condition using $\SL_3$-invariants. The expression is relatively sparse, with $110$ terms total in the standard basis, which has cardinality $42^2 = 1764$. 
	
	\begin{equation}\label{eq:seg33}
		\begin{aligned}
			\parbox[c]{0.92\textwidth}{\raggedright\hangafter=1\hangindent=2em$\  \scriptstyle \Seg_{3,3} = [123] [456] [789] ( 3 \gen{123} \gen{457} \gen{689} 
				- \gen{123} \gen{467} \gen{589} + 3 \gen{124} \gen{356} \gen{789} 		  - 3 \gen{124} \gen{357} \gen{689} + \gen{124} \gen{367} \gen{589} - \gen{124} \gen{368} \gen{579} - \gen{125} \gen{346} \gen{789} 		  + \gen{125} \gen{347} \gen{689} + \gen{127} \gen{348} \gen{569} 
				- \gen{134} \gen{258} \gen{679} - \gen{135} \gen{247} \gen{689} 		  + \gen{145} \gen{267} \gen{389} + \gen{147} \gen{258} \gen{369} ) 		 + [123] [457] [689] ( -3 \gen{123} \gen{456} \gen{789} 
				+ \gen{124} \gen{368} \gen{579} - \gen{126} \gen{348} \gen{579} 		  + \gen{135} \gen{246} \gen{789} - \gen{146} \gen{258} \gen{379} ) 		 + [123] [458] [679] ( -\gen{124} \gen{367} \gen{589} 
				- \gen{125} \gen{346} \gen{789} + \gen{126} \gen{347} \gen{589} 		  + \gen{146} \gen{257} \gen{389} ) 		 + [123] [467] [589] ( \gen{123} \gen{456} \gen{789} 
				- \gen{124} \gen{358} \gen{679} + \gen{125} \gen{348} \gen{679} 		  + \gen{134} \gen{256} \gen{789} - \gen{135} \gen{246} \gen{789} 
				+ \gen{145} \gen{268} \gen{379} ) 		 + [124] [356] [789] ( -3 \gen{123} \gen{456} \gen{789} 
				+ \gen{123} \gen{468} \gen{579} + \gen{135} \gen{247} \gen{689} 		  - \gen{135} \gen{267} \gen{489} - \gen{137} \gen{258} \gen{469} ) 		 + [124] [357] [689] ( 3 \gen{123} \gen{456} \gen{789} 
				- \gen{123} \gen{468} \gen{579} - \gen{135} \gen{246} \gen{789} 		  + \gen{136} \gen{258} \gen{479} ) 		 + [124][358][679] 
				(\gen{123}\gen{467}\gen{589} - \gen{136}\gen{257}\gen{489} ) 		+ [124][367][589] 
				(- \gen{123}\gen{456}\gen{789} + \gen{123}\gen{458}\gen{679} 		  + \gen{135}\gen{246}\gen{789} - \gen{135}\gen{268}\gen{479} ) 		+ [124][368][579] 
				(\gen{123}\gen{456}\gen{789} - \gen{123}\gen{457}\gen{689} 		  + \gen{135}\gen{267}\gen{489} ) 		+ [125][346][789] 
				(\gen{123}\gen{456}\gen{789} + \gen{123}\gen{458}\gen{679} 		  - \gen{123}\gen{468}\gen{579} - \gen{134}\gen{257}\gen{689} 		  + \gen{134}\gen{267}\gen{589} + \gen{137}\gen{248}\gen{569} ) 		+ [125][347][689] 
				(- \gen{123}\gen{456}\gen{789} + \gen{123}\gen{468}\gen{579} 		  + \gen{134}\gen{256}\gen{789} - \gen{136}\gen{248}\gen{579} ) 		+ [125][348][679] 
				(- \gen{123}\gen{467}\gen{589} + \gen{136}\gen{247}\gen{589} ) 		+ [125][367][489] 
				(- \gen{134}\gen{256}\gen{789} + \gen{134}\gen{268}\gen{579} ) 		- [125][368][479] \gen{134}\gen{267}\gen{589} 		+ [126][347][589] 
				(- \gen{123}\gen{458}\gen{679} + \gen{135}\gen{248}\gen{679} ) 
				+ [126][348][579] 
				(\gen{123}\gen{457}\gen{689} - \gen{135}\gen{247}\gen{689} ) 		- [126][357][489]\gen{134}\gen{258}\gen{679}+ [126][358][479] \gen{134}\gen{257}\gen{689} + [127][348][569] 
				(- \gen{123}\gen{456}\gen{789} + \gen{135}\gen{246}\gen{789})
				- [127][358][469] \gen{134}\gen{256}\gen{789} + [134][256][789] 
				(- \gen{123}\gen{467}\gen{589} - \gen{125}\gen{347}\gen{689} 		  + \gen{125}\gen{367}\gen{489} + \gen{127}\gen{358}\gen{469} ) 		+ [134][257][689] 
				(\gen{125}\gen{346}\gen{789} - \gen{126}\gen{358}\gen{479} ) 		+ [134][258][679] 
				(\gen{123}\gen{456}\gen{789} + \gen{126}\gen{357}\gen{489} ) 		+ [134][267][589] 
				(- \gen{125}\gen{346}\gen{789} + \gen{125}\gen{368}\gen{479} ) 	-[134][268][579] \gen{125}\gen{367}\gen{489} + [135][246][789] 
				(- \gen{123}\gen{457}\gen{689} + \gen{123}\gen{467}\gen{589} 		  + \gen{124}\gen{357}\gen{689} - \gen{124}\gen{367}\gen{589} 		  - \gen{127}\gen{348}\gen{569} ) 		+ [135][247][689] 
				(\gen{123}\gen{456}\gen{789} - \gen{124}\gen{356}\gen{789} 		  + \gen{126}\gen{348}\gen{579} ) - [135][248][679] \gen{126}\gen{347}\gen{589} + [135][267][489] 
				(\gen{124}\gen{356}\gen{789} - \gen{124}\gen{368}\gen{579} ) 		+ [135][268][479] \gen{124}\gen{367}\gen{589} - [136][247][589]\gen{125}\gen{348}\gen{679} 		+ [136][248][579] \gen{125}\gen{347}\gen{689} + [136][257][489]\gen{124}\gen{358}\gen{679}  		- [136][258][479]\gen{124}\gen{357}\gen{689} - [137][248][569] 
				\gen{125}\gen{346}\gen{789} + [137][258][469] \gen{124}\gen{356}\gen{789} 		+ [145][267][389] 
				(- \gen{123}\gen{456}\gen{789} + \gen{123}\gen{468}\gen{579} ) 	-[145][268][379] \gen{123}\gen{467}\gen{589} - [146][257][389] \gen{123}\gen{458}\gen{679}  		+ [146][258][379] \gen{123}\gen{457}\gen{689}-[147][258][369] \gen{123}\gen{456}\gen{789} .$}
		\end{aligned}
	\end{equation}
	
	\begin{remark}
		In principle, one may also compute the rank condition on $F$ in terms of brackets. However, the degree in brackets gets quite large. For eight points in $\bP^2,$ by Cramer's rule each entry of $F$ is a $8 \times 8$ minor of the Segre matrix of the eight points. Thus one may express the rank condition as a polynomial of bi-degree $(8,8)$ in the brackets $[I]$ and $\langle J \rangle$ on two copies of $\mathcal{B}_{3,8}.$ 
	\end{remark}
	
	We close this section with a theorem which states that in general, existence of a common recovery implies that the fundamental matrix has rank at most two. When there are $k^2$ or more points, there is an additional condition that the Segre matrix must drop rank. Theorem \ref{thm:vision} is a consequence of the more general Theorem 2 in \cite{bertolini} (thanks to Timothy Duff for pointing this out). However, we include a self-contained proof for our setting.
	
	\begin{theorem}\label{thm:vision}
		Let $A$ and $B$ be distinct configurations of $k^2-1$ ordered points in $\bP^{k-1}.$ Suppose that there exists a configuration of ordered points $C$ in $\bP^{2k-3}$ and linear projections $\pi_1, \pi_2: \bP^{2k-3} \dashrightarrow \bP^{k-1}$ such that $\pi_1(C) = A$ and $\pi_2(C) = B.$ Furthermore, suppose that the points are sufficiently generic for the kernel of the $k^2 \times (k^2-1)$ Segre matrix to be one-dimensional. Then the kernel $F$, when viewed as a $k \times k$ matrix, has rank at most two.
	\end{theorem}
	\begin{proof}
		Our strategy is to construct a matrix $F$ of rank at most two, such that for each pair $(a,b)$ of points in $\bP^{k-1}$ which are projections of a common point in $\bP^{2k-3},$ the product $b^TFa$ vanishes. Since $b^TFa = \sum_{i,j} F_{ij}a_ib_j,$ this would imply that the $1 \times k^2$ flattening of $F$ is in the left kernel of the $k^2 \times (k^2-1)$ Segre matrix. We first fix a partial inverse $\phi$ to $\pi_1.$ Define $V := \widehat{\pi_2(\ker \pi_1)},$ which is a $(k-2)$-dimensional subspace of $\bR^k.$ Let $M$ be the map $\pi_2 \circ \phi.$ 
		In summary, we have the data in the following commutative diagram:
		\begin{center}
			\begin{tikzcd}[column sep=small]
				& \bP^{2k-3} \arrow[dl, dashrightarrow, "\pi_1"] \arrow[dr, dashrightarrow, "\pi_2"] & \\
				\bP^{k-1} \arrow[ur, bend left, "\phi"] \arrow[rr, dashrightarrow, "M"] & & \bP^{k-1}
			\end{tikzcd}
		\end{center}
		Then consider the map
		\begin{align*}
			\tilde{F}: Gr(1,k) & \dashrightarrow Gr(k-1,k) \\
			x & \mapsto V \oplus M(x).
		\end{align*}
		This map is linear in $x$ in the sense that it is component-wise given by linear forms. More precisely, representing $V$ with a $(k-2) \times k$ matrix and $M(x)$ with a $1 \times k$ matrix of the same names, the dual Pl\"ucker coordinates of $\tilde{F}(x)$ are the $(k-1)$-subdeterminants of the $(k-1) \times k$ matrix
		\begin{equation}
			\left[
			\begin{array}{ccc}
				& M(x) &\\
				\hline
				& &
				\\			& V & 
			\end{array}
			\right]
		\end{equation}
		
		The map $\tilde{F}$ is thus represented by a matrix, which we will call $F.$ Its base locus is $\bP(V).$ Thus the rank of $F$ is two and the image is all hyperplanes in $\bP^{k-1}$ containing $\bP(V) = \bP(\ker F).$ 
		Now, suppose that there exists $c$ such that $\pi_1(c) = a$ and $\pi_2(c) = b.$ Then $c$ is in the projectivization of the vector space $\phi(a) \oplus \ker \pi_1.$ So $\pi_2(c)=b$ is on the line $F(a).$ In terms of matrix representatives, we then have $b^T Fa = 0.$
	\end{proof}

	\begin{remark}
		Theorem \ref{thm:vision} gives a necessary condition for two point configurations to be projections of a common configuration. Theorem 2 of \cite{bertolini} proves this is also sufficient, and Section 5.1.2 gives an algorithm for reconstructing the points and  projection matrices.
	\end{remark}
	
	\section{The Chow-Lam Form} \label{sec:chowlam}
	The Grassmannian $\Gr(k,n)$ is a smooth projective variety parameterizing $k$-dimensional subspaces of an $n$-dimensional vector space. Each $k \times n$ matrix $A$ determines a point of $\Gr(k,n)$ via its rowspan. The Grassmannian is embedded into projective space by taking the maximal minors of matrix representatives; these are known as \emph{(dual) Pl\"ucker coordinates} of the point. 
	
	The rest of this section will explain the Chow-Lam form for subvarieties of the Grassmannian, building towards Theorem \ref{thm:clsegre}. The Chow-Lam form was introduced in \cite{ChowLam} as a generalization of the Chow form from classical algebraic geometry. We briefly recall the setup here. 
	\begin{figure}[!h]
		\begin{center}
			\begin{tikzcd}
				& \text{Fl}(k, n-r+k,n) \arrow{dl}[swap, dash]{\pi_1} \arrow{dr}[dash]{\pi_2} & \\
				\cV \subset \Gr(k,n) & & \Gr(n-r+k,n) 
			\end{tikzcd}
		\end{center}
		\caption{Definition of the Chow-Lam locus}\label{fig:cldef}
	\end{figure}
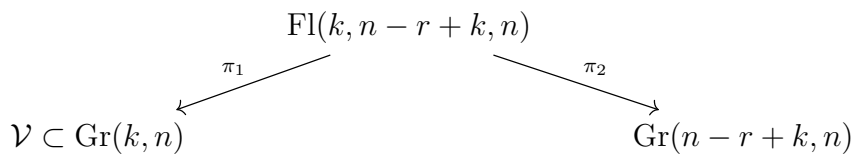
	The input data is a projective variety $\cV \subset \Gr(k,n)$ of dimension $k(r-k)-1$ for some $r \leq n.$ Then one obtains the projections in Figure \ref{fig:cldef}, where $\text{Fl}(k, n-r+k,n)$ is the partial flag variety parameterizing pairs of a $k$-space contained in an $(n-r+k)$-space in $\bC^n$. The \emph{Chow-Lam locus} $\mathcal{CL}_\cV$ is defined as the Zariski closure of $\pi_2(\pi_1^{-1}(\cV))$ in $\Gr(n-r+k,n)$. By Lemma 3.2 of \cite{ChowLam}, if $\cV$ is irreducible then $\mathcal{CL}_\cV$ is a proper irreducible subvariety of $\Gr(n-r+k, n).$ When $\mathcal{CL}_\cV$ is a hypersurface, it is defined by a single equation in Pl\"ucker coordinates, which we call the \emph{Chow-Lam form} and denote $CL_\cV$. If $\mathcal{CL}_\cV$ is not a hypersurface, we set $CL_\cV := 1.$
	
	\begin{example}[Curve in $\Gr(2,4)$]\label{eg:gr24}
		Let $\cV$ be a curve in the Grassmannian $\Gr(2,4).$  Here $k = 2, \, n = 4, $ and $r = 3.$ Thus the Chow-Lam locus lives in $\Gr(3,4) = (\bP^3)^\vee.$ It consists of planes $P$ containing a line $Q$ in $\bP^3$, such that $Q$ is a point of $\cV$ in $\Gr(2,4).$ 
		
		\begin{figure}[!h]
			\begin{minipage}{0.5\textwidth}
				\centering
					\tikzset{every picture/.style={line width=0.75pt}} 
	
	\begin{tikzpicture}[x=0.75pt,y=0.75pt,yscale=-1,xscale=1]
		
		\draw   (249,63) .. controls (249,51.95) and (273.85,43) .. (304.5,43) .. controls (335.15,43) and (360,51.95) .. (360,63) .. controls (360,74.05) and (335.15,83) .. (304.5,83) .. controls (273.85,83) and (249,74.05) .. (249,63) -- cycle ;
		\draw   (249,203) .. controls (249,191.95) and (273.85,183) .. (304.5,183) .. controls (335.15,183) and (360,191.95) .. (360,203) .. controls (360,214.05) and (335.15,223) .. (304.5,223) .. controls (273.85,223) and (249,214.05) .. (249,203) -- cycle ;
		\draw    (249,63) .. controls (262,144.5) and (262,125.5) .. (249,203) ;
		\draw    (360,63) .. controls (342,145.5) and (341,126.5) .. (360,203) ;
		\draw [color={rgb, 255:red, 74; green, 144; blue, 226 }  ,draw opacity=1 ]   (320.25,82.25) -- (249,203) ;
		\draw [color={rgb, 255:red, 208; green, 2; blue, 27 }  ,draw opacity=1 ]   (354.25,71.75) -- (270.25,218.25) ;
		\draw [color={rgb, 255:red, 74; green, 144; blue, 226 }  ,draw opacity=1 ]   (345.75,136.75) -- (300,223) ;
		\draw [color={rgb, 255:red, 74; green, 144; blue, 226 }  ,draw opacity=1 ]   (289.75,82.75) -- (257,154.5) ;
		\draw  [color={rgb, 255:red, 126; green, 211; blue, 33 }  ,draw opacity=1 ][fill={rgb, 255:red, 184; green, 233; blue, 134 }  ,fill opacity=0.2 ] (221,50.5) -- (381.5,50.5) -- (381.5,238.5) -- (221,238.5) -- cycle ;
		
		\draw (363.83,77.83) node [anchor=north west][inner sep=0.75pt]   [align=left] {$\displaystyle \textcolor[rgb]{0.49,0.83,0.13}{P}$};
		\draw (284.67,141) node [anchor=north west][inner sep=0.75pt]   [align=left] {$\displaystyle \textcolor[rgb]{0.82,0.01,0.11}{Q}$};
		\draw (326.33,193.17) node [anchor=north west][inner sep=0.75pt]   [align=left] {$\displaystyle \textcolor[rgb]{0.29,0.56,0.89}{X_{\mathcal{V}}}$};

	\end{tikzpicture}
	
			\end{minipage}
			\begin{minipage}{0.45\textwidth}
				\centering
					\tikzset{every picture/.style={line width=0.75pt}} 
	
	\begin{tikzpicture}[x=0.75pt,y=0.75pt,yscale=-1,xscale=1]
		
		\draw [color={rgb, 255:red, 126; green, 211; blue, 33 }  ,draw opacity=1 ]   (220.48,79.69) -- (398.19,200.98) ;
		\draw [shift={(400.67,202.67)}, rotate = 214.31] [fill={rgb, 255:red, 126; green, 211; blue, 33 }  ,fill opacity=1 ][line width=0.08]  [draw opacity=0] (8.93,-4.29) -- (0,0) -- (8.93,4.29) -- cycle    ;
		\draw [shift={(218,78)}, rotate = 34.31] [fill={rgb, 255:red, 126; green, 211; blue, 33 }  ,fill opacity=1 ][line width=0.08]  [draw opacity=0] (8.93,-4.29) -- (0,0) -- (8.93,4.29) -- cycle    ;
		\draw [color={rgb, 255:red, 74; green, 144; blue, 226 }  ,draw opacity=1 ]   (209.55,216.79) .. controls (227.98,69.82) and (276.35,224.67) .. (311,142) ;
		\draw [shift={(209,221.33)}, rotate = 276.64] [fill={rgb, 255:red, 74; green, 144; blue, 226 }  ,fill opacity=1 ][line width=0.08]  [draw opacity=0] (8.93,-4.29) -- (0,0) -- (8.93,4.29) -- cycle    ;
		\draw [color={rgb, 255:red, 208; green, 2; blue, 27 }  ,draw opacity=1 ]   (338.6,63) -- (311,142) ;
		\draw [shift={(311,142)}, rotate = 109.26] [color={rgb, 255:red, 208; green, 2; blue, 27 }  ,draw opacity=1 ][fill={rgb, 255:red, 208; green, 2; blue, 27 }  ,fill opacity=1 ][line width=0.75]      (0, 0) circle [x radius= 3.35, y radius= 3.35]   ;
		\draw [color={rgb, 255:red, 74; green, 144; blue, 226 }  ,draw opacity=1 ]   (338.8,63.03) -- (311,142) ;
		\draw [shift={(339.8,60.2)}, rotate = 109.4] [fill={rgb, 255:red, 74; green, 144; blue, 226 }  ,fill opacity=1 ][line width=0.08]  [draw opacity=0] (8.93,-4.29) -- (0,0) -- (8.93,4.29) -- cycle    ;
		
		\draw (317.67,126.33) node [anchor=north west][inner sep=0.75pt]   [align=left] {$\displaystyle \textcolor[rgb]{0.82,0.01,0.11}{Q}$};
		\draw (347.33,207.67) node [anchor=north west][inner sep=0.75pt]   [align=left] {$\displaystyle \textcolor[rgb]{0.49,0.83,0.13}{\Gr( 2,P)}$};
		\draw (221,185.33) node [anchor=north west][inner sep=0.75pt]   [align=left] {$\displaystyle \mathcal{\textcolor[rgb]{0.29,0.56,0.89}{V}}$};

	\end{tikzpicture}
	
			\end{minipage}
			\caption{Geometry in $\bP^3$ (left) and $\Gr(2,4)$ (right)}
			\label{fig:twoschubert}
		\end{figure}
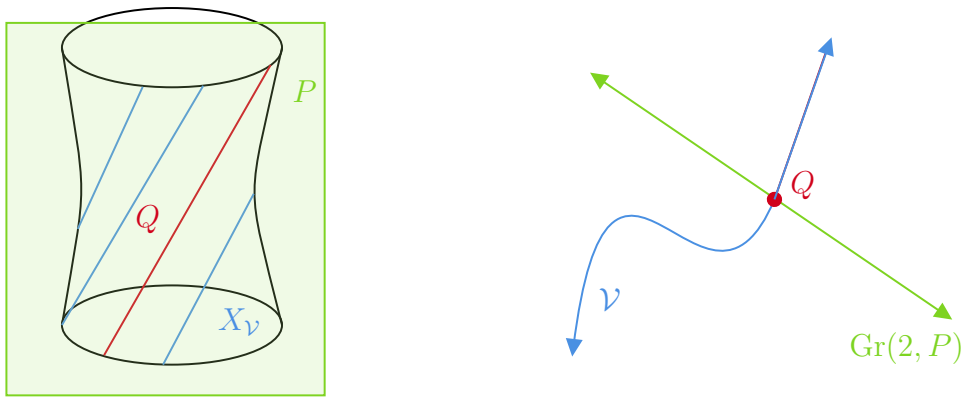
		
		Let $X_\cV$ be the ruled surface in $\bP^3$ swept out by all of the lines in $\cV$. Then any tangent plane to $X_\cV$ contains the two lines in $X_\cV$ through the point of tangency; thus the dual surface $X_\cV^\vee$ is contained in $\mathcal{CL}_\cV.$ In fact, since they are irreducible varieties of the same dimension, we have that the Chow-Lam locus is precisely $X_\cV^\vee.$
	\end{example}
	
	The case $k=1$ recovers the classical Chow form for projective varieties. It was introduced by Chow and van der Waerden as a way to encode any projective variety by a single polynomial \cite{CvW}. That is, given a Chow form of $\cV,$ one may recover the radical ideal of $\cV.$ They used this to construct a space known as the \emph{Chow variety} parameterizing cycles of a fixed degree and dimension in a projective space. For an introduction to Chow forms see \cite{dalbec}.
	
	Unique recovery of a variety fails for the Chow-Lam form even for curves in $\Gr(2,4)$, as one can see in Example \ref{eg:gr24}. The two rulings of $X_\cV$ give two different curves in $\Gr(2,4)$ with the same Chow-Lam locus. The extent to which one may recover a subvariety of $\Gr(k,n)$ from its Chow-Lam form is explored in \cite{PrattRanestad}. In particular, for curves whose corresponding ruled surface is not a cone, it turns out that ruling the same surface is equivalent to having the same Chow-Lam form.
	
	To compute Chow-Lam forms, it is convenient to use multiple different coordinate systems on the Grassmannian. One may parameterize a $k$-subspace $A$ of $\bC^n$ as the rowspan of a matrix. Then the entries of this matrix are called the \emph{dual Stiefel coordinates} and the maximal minors are the \emph{dual Pl\"ucker coordinates}. Alternatively, one may represent $A$ as the kernel of an $(n-k) \times n$ matrix.  In this case the entries are called the \emph{primal Stiefel coordinates} and the maximal minors are called the \emph{primal Pl\"ucker coordinates} of $A.$ 
	
	We denote the dual Pl\"ucker coordinates by $q_I(A)$ and the primal Pl\"ucker coordinates by $p_J(A),$ where $I \subset \binom{[n]}{k}$ and $J \subset \binom{[n]}{n-k}.$ They are related up to sign by taking complements of indices: specifically, if $I = [i_1 \cdots i_k]$ we define \[\varepsilon(I) := i_1 + \ldots + i_k - (1 + \ldots + k).\] Then $p_I = (-1)^{\varepsilon(I)}q_{[n] \setminus I}.$  If $k$ is close to $n$ then it is more convenient to use primal coordinates. For example, the ten Pl\"ucker coordinates on $\Gr(3, 5)$ are
	
	\begin{equation}
		\label{eq:primaldual35}
		\begin{matrix}
			p_{12} & p_{13} & p_{14} & p_{15} & p_{23} & p_{24} & p_{25} & p_{34} &
			p_{35} & p_{45}, \\ 
			q_{345} & -q_{245} & q_{235} & -q_{234} & q_{145} & -q_{135} & q_{134} & q_{125} & -q_{124} & q_{123}.
		\end{matrix}
	\end{equation}
	
	An important statistic of a subvariety of a Grassmannian is its class in the cohomology ring $H^*(\Gr(k,n), \bZ).$ As an abelian group, this is isomorphic to $\bZ^{\binom{n}{k}}$ with basis given by the \emph{Schubert cycles}. We will define these with respect to intersections with the standard flag and index them using partitions. 
	
	Let $e_1, \, \ldots, \, e_n$ be the standard basis of $\bC^n,$ and let $E_i := \text{span}(e_1, \ldots, e_i).$ For a partition $\lambda = (\lambda_1, \ldots, \lambda_k)$ fitting inside a $k \times (n-k)$ box, define the Schubert variety $\Omega_\lambda$ to be
	\begin{equation}\label{eqn:schubert}
		\Omega_\lambda = \{L \in \Gr(k,n) \ : \ \dim \, L \cap E_{n-k+\lambda_i - i} \geq i\}.
	\end{equation}
	The Schubert variety $\Omega_\lambda$ is a closed irreducible subvariety of codimension $\sum_i \lambda_i.$ In particular,  $\Omega_\varnothing$ is the class of $\Gr(k,n)$ and $\Omega_{k \times (n-k)}$ is the class of a point. Every subvariety $\cV$ of $\Gr(k,n)$ has a unique expansion into Schubert classes:
	\begin{align*}
		[\cV] = \sum_{\lambda \, \subseteq  \, k \times (n-k)} \delta_\lambda (\cV) \cdot [ \Omega_\lambda].
	\end{align*} 
	Let $\cV$ be a subvariety of $\Gr(k,n)$ of dimension $k(r-k)-1.$ Let $\alpha = (n-r+1, n-r, \ldots , n-r).$ This is the unique partition $\alpha$ with size equal to the codimension of $\cV$ satisfying $\alpha_1 = n-r+1.$ We set $\alpha(\cV) := \delta_\alpha(\cV)$ and call this the \emph{Chow-Lam degree} of~$\cV.$ 
	
	\begin{lemma}[Theorem 3.5 of \cite{ChowLam}]
		Let $\cV$ be a subvariety of dimension $k(r-k)-1$ in the Grassmannian ${\rm Gr}(k,n)$. The Chow-Lam form $CL_\cV$ is a polynomial of degree $\alpha(\cV)$ in the Pl\"ucker coordinates on ${\rm Gr}(k+n-r, n)$. 
	\end{lemma} 
	
	\section{Torus Orbits in the Grassmannian}\label{sec:torus}
	
	In this section we study torus orbits in the Grassmannian and their Chow-Lam forms. The Grassmannian $\Gr(k,n)$ is equipped with an action of $T := (\bC^*)^n.$ In terms of a matrix parameterization, this may be seen as scaling the columns of a $k \times n$ matrix representative. Given any point $A \in \Gr(k,n),$ we write $\cT_A$ for the Zariski closure of the orbit $T \cdot A$ in $\Gr(k,n).$ This is a toric variety of dimension at most $n-1.$ It has dimension exactly $n-1$ if $A$ is general, in particular if the Pl\"ucker coordinates of $A$ are all non-zero. As in the introduction, we fix positive integers $k$ and $l$ and let $n := k\ell.$
	
	\begin{example}[Torus orbit in $\Gr(2,6)$]\label{eg:cl23}
		Let $k=2$ and $\ell=3,$ and let $A$ be a generic point in $\Gr(2,6)$ whose Pl\"ucker coordinates are non-zero. Then $\cT_A$ has dimension $6-1 = 5.$ The Chow-Lam locus lives in $\Gr(3,6).$ Following \cite[Example 4.3]{ChowLam}, it is parameterized in dual Stiefel coordinates by the matrices
		$$
		\begin{bmatrix}
			a_{11} t_1 & a_{12} t_2 & a_{13} t_3 & a_{14} t_4 & a_{15} t_5 & a_{16} t_6 \\
			a_{21} t_1 & a_{22} t_2 & a_{23} t_3 & a_{24} t_4 & a_{25} t_5 & a_{26} t_6 \\
			y_1 & y_2 & y_3 & y_4 & y_5 & y_6 
		\end{bmatrix}.
		$$
		Here $(t_1, ..., t_6)$ varies over elements of $(\bC^*)^6.$ However, we could also parameterize it in primal Stiefel coordinates as $3 \times 6$ matrices $B$ such that, for some $t \in (\bC^*)^6,$ we have
		$$\begin{bmatrix}
			b_{11} & b_{12} & b_{13} & b_{14} & b_{15} & b_{16} \\
			b_{21} & b_{22} & b_{23} & b_{24} & b_{25} & b_{26}\\
			b_{31} & b_{32} & b_{33} & b_{34} & b_{35} & b_{36}
		\end{bmatrix} \cdot \text{diag}(t_1, ..., t_6) \cdot  \begin{bmatrix}
			a_{11} & a_{21} \\
			a_{12} & a_{22} \\
			a_{13} & a_{23} \\
			a_{14} & a_{24} \\
			a_{15} & a_{25} \\
			a_{16} & a_{26}
		\end{bmatrix} = 0.$$
		Re-arranging, we obtain the expression
		\[t_1 \begin{bmatrix}
			a_{11}b_{11} \\
			a_{11}b_{21} \\
			a_{11}b_{31} \\
			a_{21}b_{11} \\
			a_{21}b_{21} \\
			a_{21}b_{31}
		\end{bmatrix} +t_2 \begin{bmatrix}
			a_{12}b_{12} \\
			a_{12}b_{22} \\
			a_{12}b_{32} \\
			a_{22}b_{12} \\
			a_{22}b_{22} \\
			a_{22}b_{32}
		\end{bmatrix} + t_3\begin{bmatrix}
			a_{13}b_{13} \\
			a_{13}b_{23} \\
			a_{13}b_{33} \\
			a_{23}b_{13} \\
			a_{23}b_{23} \\
			a_{23}b_{33}
		\end{bmatrix} + t_4\begin{bmatrix}
			a_{14}b_{14} \\
			a_{14}b_{24} \\
			a_{14}b_{34} \\
			a_{24}b_{14} \\
			a_{24}b_{24} \\
			a_{24}b_{34}
		\end{bmatrix} + t_5\begin{bmatrix}
			a_{15}b_{15} \\
			a_{15}b_{25} \\
			a_{15}b_{35} \\
			a_{25}b_{15} \\
			a_{25}b_{25} \\
			a_{25}b_{35}
		\end{bmatrix} + t_6 \begin{bmatrix}
			a_{16}b_{16} \\
			a_{16}b_{26} \\
			a_{16}b_{36} \\
			a_{26}b_{16} \\
			a_{26}b_{26} \\
			a_{26}b_{36}
		\end{bmatrix} = 0.\]
		Thus $B$ is in the Chow-Lam locus of $\cT_A$ in $\Gr(3,6)$ if and only if the Segre determinant $\text{Seg}_{2,3}(A,B)$ vanishes. In dual Pl\"ucker coordinates on $\Gr(2,6)$ and primal Pl\"ucker coordinates on $\Gr(3,6),$ we obtain the expression 
		\[\begin{matrix} \Seg_{2,3} & = &
			([12][34][56] \,+\, [14][25][36])\, \gen{123} \gen{456}
			\,-\,  [13] [25] [46] \, \gen{124} \gen{356} \\ & & 
			\,+\,  [12] [35] [46] \, \gen{134} \gen{256}
			\,-\,  [12] [34] [56] \, \gen{135} \gen{246}
			\,+\,  [13] [24] [56] \, \gen{125} \gen{346}. \qedhere
		\end{matrix} \]
	\end{example}
	
	Suppose that $n=k\ell$ for some $\ell \geq 2.$ Let $A$ be a general point in $\Gr(k,n)$. Then the Chow-Lam locus is a subvariety of $\Gr(n-\ell,n).$ The analysis in Example \ref{eg:cl23} extends to the following result.
	
	\begin{theorem}[Segre Determinant]\label{thm:clsegre}
		Suppose $k, \ell \geq 2$ and let $n = k\ell.$ Fix a point $A$ in ${\rm Gr}(k,n)$ with non-zero Pl\"ucker coordinates. Then the Chow-Lam form of $\cT_A$ in primal Pl\"ucker coordinates $B$ on ${\rm Gr}(n-\ell, n)$ equals the Segre determinant ${\rm Seg}_{k,\ell}(A,B).$
	\end{theorem}
	
	We will need tools from later in this section to prove Theorem \ref{thm:clsegre} completely. For now, we establish the following lemma.
	
	\begin{lemma}[Factor of Segre Determinant]\label{lem:cltorus}
		Fix $k, \ell \geq 2$ and let $n = k\ell.$ Fix a point $A$ in ${\rm Gr}(k,n)$ such that $\dim \cT_A = n-1.$ Then the Chow-Lam form of $\cT_A$ in primal Pl\"ucker coordinates $B$ on ${\rm Gr}(n-\ell, n)$ divides the Segre determinant ${\rm Seg}_{k,\ell}(A,B).$ 
	\end{lemma}
	
	\begin{proof}
		The Chow-Lam locus is the Zariski closure of the set of points $B$ in $\Gr(n-\ell,n)$ which contain a subspace $t \cdot A$ for some $t \in T.$ Representing $B$ in primal Stiefel coordinates and $t \cdot A$ in dual Stiefel coordinates, we are seeking $\ell \times n$ matrices $B$ such that $B \cdot \text{diag}(t_1, \ \ldots \ , t_n) \cdot A = 0.$ Re-arranging, we obtain the condition that for some $t \in (\bC^*)^n,$
		\begin{equation}\label{eq:dependence}
			\sum_{i = 1}^n t_i (A_i \otimes B_i) = 0.
		\end{equation}
		If \eqref{eq:dependence} holds, then $A_1 \otimes B_1, \ \ldots, \ A_n \otimes B_n$ are linearly dependent. Thus the Chow-Lam form is an irreducible factor of the Segre determinant. 
	\end{proof}
	
	\begin{example} \label{eg:clfactors}
		For special matrices $A,$ the Segre determinant becomes reducible and the Chow-Lam form is one of the irreducible factors. For instance, let $A$ in $\Gr(2,4)$ be any point whose Pl\"ucker coordinate $[12]$ is zero. Then $\cT_A$ is the $3$-dimensional variety with ideal generated by $q_{12}$ and the Pl\"ucker relation for $\Gr(2,4)$. By Example \ref{eg:22}, 
		\begin{align*}
			\Seg_{2,2}(A, B) = [13][24]\gen{12}\gen{34}.
		\end{align*}
		However, the Chow-Lam form of $\cT_A$ in primal coordinates is $\gen{34}.$ Indeed, the vanishing of $p_{34} = q_{12}$ exactly cuts out the original variety. The extra factor of $\gen{12}$ represents the lines passing through the singular point $A' := \text{span}(e_3, e_4)$ in $\cT_A.$ Indeed, any matrix with $p_{34} = 0$ satisfies $B \cdot A' = 0.$ In terms of the proof of Lemma \ref{lem:cltorus}, we have a linear dependence where the coefficients $(0,0,t_3,t_4)$ are not in $(\bC^*)^4.$
	\end{example}
	
	To better predict situations like Example \ref{eg:clfactors}, we will introduce some combinatorial tools to compute the Chow-Lam degree of $\cT_A.$ The \emph{matroid} of a point $A$ in the Grassmannian $\Gr(k,n)$ is defined by its \emph{bases}, namely the collection of indices $I \in \binom{[n]}{k}$ such that the corresponding Pl\"ucker coordinate $q_I$ is nonzero. In general, the variety $\cT_A$ will (up to isomorphism) only depend on the underlying matroid of the point $A,$ see e.g. \cite[Proposition 13.12]{Nonlinear}. Thus we may denote $\cT_A$ by $\cT_M,$ where $M$ is the matroid of $A.$ We let $\delta_\lambda(M)$ denote the Schubert coefficient $\delta_\lambda(\cT_M)$ and call the collection of these the \emph{Schubert coefficients of the matroid $M$}. In particular, $\alpha(M)$ denotes the Chow-Lam degree of $\cT_M.$ 
	
	The \emph{uniform matroid} $U_{k,n}$ has as its bases all size $k$ subsets of $[n].$ It arises as the matroid of a point whose Pl\"ucker coordinates are all nonzero. Klyachko \cite{klyachko} computed the Schubert coefficients of the uniform matroid in terms of dimensions of irreducible $SL_n$-representations. His formula involves a count $\# \text{SSYT}(\lambda, n)$ of the number of semi-standard Young tableaux of shape $\lambda$ with entries in $[n].$ This formula comes from representation theory; in that context, the number $\text{SSYT}(\lambda, n)$ is the dimension of the irreducible $SL_n$-representation obtained by applying the Schur functor $\bS_\lambda$ to the standard representation of $SL_n.$ The partition complement $\lambda^c$ is obtained by removing $\lambda$ from a $k \times (n-k)$ rectangle and rotating by $180$ degrees.
	
	\begin{proposition}[Theorem 6 of \cite{klyachko}]\label{prop:klyachko}
		Let $\lambda$ be a partition fitting in a $k \times (n-k)$ rectangle. Then the coefficient $\delta_\lambda(U_{k,n})$ is
		\begin{equation}\label{eq:klyachko}
			\delta_\lambda(U_{k,n}) \ = \ \sum_{i = 0}^k (-1)^i \binom{n}{i} \# \text{SSYT}(\lambda^c, k-i).
		\end{equation}
	\end{proposition}
	
	\begin{proof}[Proof of Theorem~\ref{thm:clsegre}]
		By Lemma \ref{lem:cltorus}, the Chow-Lam form divides the Segre determinant. Thus it suffices to show that the Chow-Lam degree of $U_{k,n}$ is $k.$ We do this using Klyachko's formula. The complement of $\alpha = (n-r+1, n-r, \, \ldots, \, n-r)$ in the $k \times (n-k)$ rectangle is $\alpha^c =(r-k, \, \ldots, \, r-k, r-k -1).$ Because $\alpha^c$ has $k$ parts, the contribution to the sum in \eqref{eq:klyachko} is nonzero only when $i = 0.$ The semi-standard condition fixes all but the last column of each tableau of shape $\alpha^c$, giving that $\alpha(U_{k,n}) = k.$ 
	\end{proof}
	
	\section{The Segre Coefficient Variety}\label{sec:chowquotient}
	In this section, we introduce the \emph{Segre coefficient variety}, which parameterizes Segre determinants as polynomials in the $B$ variables as $A$ varies. We prove in Theorem \ref{thm:independent} that the linear span of the monomials appearing in the Segre coefficient map is as large as possible. This results in Corollary \ref{cor:git}, which states that for $k = 2,$ the Segre coefficient variety recovers the GIT quotient $(\bP^1)^{2 \ell} / \!/ \text{SL}(2)$ parameterizing configurations of points in $\bP^1$. 

	Let $\Gr(k,k\ell)^\circ \subset \Gr(k,k\ell)$ be the Zariski open subset of points whose matroid is uniform. Consider the map
	\begin{equation}\label{eq:segcoeffmap}
		\begin{split}
			\pi: \Gr(k,k\ell)^\circ & \to \bP H^0(\Gr(\ell, k\ell), \cO_{\Gr(\ell,k\ell)}(k)) \\
			A & \mapsto \Seg(A,B).
		\end{split}
	\end{equation}
	The map $\pi$ sends a point $A$ to $\Seg(A,B),$ viewed as a polynomial in the $B$-coordinates. While $\pi$ is defined without choosing a basis for the target space, it is often convenient to take a basis of standard $B$-monomials of degree $k$ for $\cO_{\Gr(\ell,k\ell)}(k)$, as in Example \ref{eg:segcubic}. Then $\pi$ sends $A$ to the vector of $A$-coefficients of the standard $B$-monomials appearing in $\Seg_{k,\ell}(A,B)$. We define the $\emph{Segre coefficient variety}$ $\text{Coeff}(\Seg_{k,\ell})$ as the Zariski closure of the image of $\pi$.  
	\begin{example}[Segre cubic]\label{eg:segcubic}
		From Example \ref{eg:cl23} we get the map
		\begin{equation}
			\begin{split}
				\pi: \Gr(2,6)^\circ & \to \bP^4 \\
				A & \mapsto ([12][34][56] +[14][25][36], \, \\
				& -[13] [25] [46], \ [12] [35] [46], \ -[12] [34] [56], \ [13] [24] [56]).	
			\end{split}
		\end{equation}
		Its image is cut out by the following degree $3$ polynomial, which is known as the \emph{Segre cubic}: 
		\begin{equation*}
			x_0x_1x_3-x_1x_2x_3-x_0x_2x_4-x_1x_2x_4-x_1x_3x_4-
			x_2x_3x_4.
		\end{equation*}
		The corresponding variety is known as the \emph{Segre cubic threefold}. It is the unique (up to an isomorphism) cubic hypersurface in $\bP^4$ with the maximum number of ordinary double points, namely ten \cite{kalker}.
	\end{example}
	
	The Segre coefficient variety recovers a known construction for $k=2$, namely the GIT quotient $(\bP^1)^{2 \ell} / \!/ \text{SL}_2$ parameterizing configurations of $2 \ell$ distinct points on a projective line. This is defined as ${\rm Proj} \, R,$ where $R$ is the graded ring of $SL_2$-invariants of $2 \ell$ ordered points in $\bP^1.$ This ring is studied by Howard, Millson, Snowden, and Vakil in \cite{invariantn}, and they give generators for the ideal of relations between these invariants. The equivalence of that construction with our construction rests on Theorem \ref{thm:independent}.

	\begin{theorem}\label{thm:independent}
		Every monomial of the form $[I_1] \cdots [I_\ell]$ with $I_1 \cup \, \ldots \, \cup I_\ell = [k\ell]$ appears in the linear span of the coefficients of the $B$-monomials in $\Seg_{k, \ell}(A,B).$
	\end{theorem}
	
	\begin{proof}
		Let $V$ be the vector space spanned by monomials $[I_1] \cdots [I_\ell]$ such that $I_1 \cup \, \ldots \, \cup I_\ell = [k\ell].$ Suppose that the coefficients of the Segre determinant do not span $V.$ Then there is a linear relation among them when $\Seg_{k, \ell}(A,B)$ is written in the basis of standard monomials. Thus the Segre coefficient variety lies in a hyperplane in $\bP(V).$ 
		
		We show that the Segre coefficient variety does not lie in a hyperplane. To do this, we will produce a collection of points in ${\rm Coeff}({\rm Seg}_{k, \ell})$ whose linear span is all of $\bP (V).$ 
		
		Consider the point $p$ in $\Gr(k,k\ell)$ given by the rowspan of a matrix where columns $(i-1) \ell + 1$ through $i\ell$ are equal to $e_i,$ for $1 \leq i \leq k.$ The Segre matrix develops a block form in this case, and the determinant (up to scalar) is $\gen{1 \cdots \ell} \gen{\ell + 1 \cdots 2\ell} \cdots \gen{(k-1)\ell + 1 \cdots k \ell}.$ Thus the Segre determinant of $\cT_p$ is nonzero, so the map $\pi$ extends from $\Gr(k,k \ell)^\circ$ to $p.$ By symmetry, we may take any point $p$ given as the rowspan of a matrix whose columns are partitioned into $k$ classes of $\ell$ pairwise parallel elements. By continuity, the image $\pi(p)$ will be in the closure of the Segre coefficient variety. The collection of these points together spans all of $\bP(V).$
	\end{proof}
	
	\begin{corollary}\label{cor:git}
		The variety ${\rm Coeff}({\rm Seg}_{2, \ell})$ is isomorphic to the GIT quotient $(\bP^1)^{2 \ell} / \!/ \rm{SL}_2.$
	\end{corollary}
	\begin{proof}
		Theorem \ref{thm:independent} tells us that every monomial $[I_1] \cdots [I_\ell]$ where $I_1 \cup \ldots \cup I_\ell = [k \ell]$ lies in the linear span of the coefficients of the Segre polynomial. The GIT quotient (with respect to the linearization $1^n$) is defined as $\rm{Proj} \, R,$ where $R$ is the ring of invariants of $\SL_2$ acting diagonally on the product of $n$ copies of $\bP^1$; see the introduction of \cite{invariantn} for a more precise description of the ring $R$. It suffices to show that these monomials generate $R$ as a ring. This is the content of Kempe's theorem \cite{Kempe}. 
	\end{proof}

	If we drop the assumption that $I_1 \cup \, \ldots \, \cup I_\ell = [k\ell]$ and instead range over all sets $\{I_1, \, \ldots, \, I_{\ell}\}$ with entries in $[k\ell],$ the collection of monomials $[I_1] \cdots [I_\ell]$ uniquely determines the torus orbit closure; see e.g. \cite[Chapter 13]{Nonlinear}. However, the map $\pi$ is not injective on generic torus orbit closures in general.
	
	\begin{example}[Six points in $\bP^2$]\label{eg:degree2}
		From Example \ref{eg:cl23} we have that the Segre coefficient map is linearly equivalent to the map
		\begin{equation}\label{eq:p2map}
			\begin{split}
				\pi': \Gr(3,6)^\circ & \to \bP^4 \\
				A & \mapsto ([123][456], \, [124][356], \, [125][345], \, [134][256], \, [135][246]).	
			\end{split}
		\end{equation}
		If we parameterize the Grassmannian by $3 \times 6$ matrices, then the following two matrices have nonzero Pl\"ucker coordinates and go to the same point $(-2,1,8,1,8)$ under the map $\pi'$:
		\[p \ = \begin{bmatrix}
			1 & 0 & 0 & 1 & 1 & 1 \\
			0 & 1 & 0 & 1 & 2 & 4 \\
			0 & 0 & 1 & 1 & 3 & 5
		\end{bmatrix}, \qquad q \ := \begin{bmatrix} 1 & 0 & 0 & 1 & 1 & 1 \\
			0 & 1 & 0 & 1 & -1 & -3 \\
			0 & 0 & 1 & 1 & -2 & -4
		\end{bmatrix}.\]
		However, their torus orbit closures are different. We may see this by noting that $[123][145][246][245]$ evaluates to $-8$ for the first point and $4$ for the second in the affine chart given by $[123] = 1$, and is constant on torus orbits.
	\end{example}
	
	Recall that by Theorem \ref{thm:clsegre}, the Segre determinant computes the Chow-Lam form of a torus orbit closure. From the point of view of the Chow-Lam form, this gives us a new family of varieties which are distinct but have the same Chow-Lam form, thus building on the work of \cite{PrattRanestad}. We hope to explore the case of $k = 3$ and general $\ell$ in future work.

	\bigskip
	\bigskip
	
	\footnotesize
	\noindent {\bf Author's address:}
	
	\smallskip
	
	\noindent Elizabeth Pratt, UC Berkeley
	\hfill \url{epratt@berkeley.edu}

\end{document}